# Regression rank scores in nonlinear models

**Jana Jurečková**[*,1]

*Charles University*

**Abstract:** Consider the nonlinear regression model

$$Y_i = g(\mathbf{x}_i, \boldsymbol{\theta}) + e_i, \quad i = 1, \ldots, n \qquad (1)$$

with $\mathbf{x}_i \in \mathbb{R}^k$, $\boldsymbol{\theta} = (\theta_0, \theta_1, \ldots, \theta_p)' \in \boldsymbol{\Theta}$ (compact in $\mathbb{R}^{p+1}$), where $g(\mathbf{x}, \boldsymbol{\theta}) = \theta_0 + \tilde{g}(\mathbf{x}, \theta_1, \ldots, \theta_p)$ is continuous, twice differentiable in $\boldsymbol{\theta}$ and monotone in components of $\boldsymbol{\theta}$. Following Gutenbrunner and Jurečková (1992) and Jurečková and Procházka (1994), we introduce regression rank scores for model (1), and prove their asymptotic properties under some regularity conditions. As an application, we propose some tests in nonlinear regression models with nuisance parameters.

## 1. Introduction

Consider the nonlinear regression model

$$(1.1) \qquad Y_i = g(\mathbf{x}_i, \boldsymbol{\theta}) + e_i, \quad i = 1, \ldots, n$$

where $\mathbf{Y} = (Y_1, \ldots, Y_n)'$ is a vector of observations, $\mathbf{x}_i \in \mathcal{X} \subset \mathbb{R}^q_+$, $i = 1, \ldots, n$ are given vectors, $e_1, \ldots, e_n$ are *i.i.d.* errors with a positive (but generally unknown) density $f$ and $\boldsymbol{\theta} = (\theta_0, \theta_1, \ldots, \theta_p)'$ is an unknown parameter. We assume that $\theta_0$ is an intercept, i.e. that $g(\mathbf{x}, \boldsymbol{\theta}) = \theta_0 + \tilde{g}(\mathbf{x}, \theta_1, \ldots, \theta_p)$. Koenker and Bassett [14] introduced the $\alpha$-regression quantile and $\alpha$-trimmed least squares estimator for the linear regression model. Their idea is very natural and regression quantiles soon became very popular among applied statisticians and econometricians. Gutenbrunner and Jurečková [3] showed that the variables, dual to regression quantiles in the parametric linear programming sense, extend the rank scores to the linear regression model; such dual regression quantiles were called the regression rank scores (RRS). They are invariant to the (linear) regression, and as such are suitable for construction of tests in the presence of nuisance regression. Such tests were considered already by Gutenbrunner and Jurečková [3], and then Gutenbrunner et al. [4] constructed a general class of tests of linear hypothesis based on regression rank scores. Koul and Saleh [17] extended the regression quantiles and regression rank scores to the linear autoregressive time series. The autoregression rank scores were then used for testing by Hallin et al. [7], Hallin et al. [8], Hallin and Jurečková [6], Kalvová et al. [12], among others.

---

[*]Supported by the Research Projects MSM 0021620839 and LC06024 and by the Czech Republic Grant 201/05/2340.

[1]Department of Probability and Statistics, Charles University, CZ-186 75 Prague 8, Czech Republic, e-mail: jurecko@karlin.mff.cuni.cz

The basic definition of the $\alpha$-regression quantile can be naturally extended to other models, including nonlinear regression. Chen [1] used the technique of the linearization of the regression function. Jurečková and Procházka [10] proved the consistency and asymptotic normality of regression quantiles in model (1.1), covering the logistic regression and a mixture of two exponentials as a regression function. There already exists a rich literature on nonlinear regression quantiles and their computation. We refer to Koenker and Park [15] and to Koenker [13], where many other references are cited.

A natural idea is to define the nonlinear regression rank scores as some form of duals to nonlinear regression quantiles. However, such dual variables do not retain the advantages of the RRS in the linear model, mainly they are not invariant to the (nonlinear) regression. Mukherjee [18], inspired by the dual steps of Koenker and Park [15] in their interior point algorithm for nonlinear regression quantiles, proposed the regression rank scores for a nonlinear time series model. His RRS are not a straightforward extension of RRS of Gutenbrunner and Jurečková [3], but under some further regularity conditions their asymptotic behavior is analogous to that of the linear RRS. Mukherjee [18] proved the asymptotic representations of the regression rank scores process and of the nonlinear rank statistics, but to a construction of tests in models affected by a nonlinear regression with nuisance parameters he missed the (even asymptotic) invariance of the regression rank scores to this type of regression.

Following Jurečková and Procházka [10], Koenker and Park [15], Mukherjee [18], and El-Attar et al. [2], we shall consider a possible version of regression rank scores in model (1.1). Our ultimate goal is their possible application in testing with nuisance (nonlinear) regression.

## 2. Regression rank scores

We shall work with the model (1.1) under the conditions of Jurečková and Procházka [10], namely:

**(A.1)** The function $g(\mathbf{x}, \boldsymbol{\theta}) : \mathcal{X} \times \boldsymbol{\theta} \mapsto \mathbb{R}^1$ has the form

$$g(\mathbf{x}, \boldsymbol{\theta}) = \theta_0 + \tilde{g}(\mathbf{x}, (\boldsymbol{\theta}^*), \ \mathbf{x} \in \mathcal{X}, \ \boldsymbol{\theta}^* = (\theta_1, \ldots, \theta_p)', \ (\theta_0, \boldsymbol{\theta}^{*\prime})' \in \boldsymbol{\Theta}$$

with some function $\tilde{g}$ of $\mathbf{x}$ and of $\boldsymbol{\theta}^*$.

We assume that function $g(\mathbf{x}, \boldsymbol{\theta})$ is strictly monotone and twice differentiable in every component of $\boldsymbol{\theta} = (\theta_0, \theta_1, \ldots, \theta_p)'$. The first and second derivatives are bounded by $K$, $0 < K < \infty$, uniformly in $\mathcal{X}$ and $\boldsymbol{\Theta}$.

**(A.2)** The parameter space $\boldsymbol{\Theta}$ and the space $\mathcal{X}$ of $\mathbf{x}$ are compact.

**(A.3)** (Identifiability). Every set $(y_1, \mathbf{x}_1), \ldots, (y_{p+1}, \mathbf{x}_{p+1})$ of $p+1$ different points determines uniquely $\boldsymbol{\theta} \in \boldsymbol{\Theta}$ such that

$$y_i = g(\mathbf{x}_i, \boldsymbol{\theta}), \quad i = 1, \ldots, p+1.$$

**(A.4)** The errors $e_1, \ldots, e_n$ are independent, identically distributed with a symmetric, positive and bounded density $f$, that has a bounded derivative $f'$.

**(A.5)** There exist finite positive constants $k_1$, $k_2$ such that, for $n > n_0$,

$$k_1 \|\boldsymbol{\theta}_2 - \boldsymbol{\theta}_1\|^2 \leq \frac{1}{n} \sum_{i=1}^{n} [g(\mathbf{x}_i, \boldsymbol{\theta}_2) - g(\mathbf{x}_i, \boldsymbol{\theta}_1)]^2 \leq k_2 \|\boldsymbol{\theta}_2 - \boldsymbol{\theta}_1\|^2$$

where $\|\cdot\|$ stands for the Euclidean norm.

**(A.6)** Put

$$v_{ij}(\boldsymbol{\theta}) = \left[\frac{\partial g(\mathbf{x}_i, \boldsymbol{\theta}+\boldsymbol{\delta})}{\partial \delta_j}\right]_{\boldsymbol{\delta}=\mathbf{0}}, \quad i=1,\ldots,n; \ j=0,1,\ldots,p$$

and denote $\mathbf{V}_n(\boldsymbol{\theta}) = \mathbf{V}(\boldsymbol{\theta}) = \left[v_{ij}(\boldsymbol{\theta})\right]_{i=1,\ldots,n}^{j=0,\ldots,p}$. We shall assume that

$$\lim_{n\to\infty} \mathbf{Q}_n(\boldsymbol{\theta}) = \mathbf{Q}(\boldsymbol{\theta})$$

where $\mathbf{Q}_n(\boldsymbol{\theta}) = \frac{1}{n}\mathbf{V}'_n(\boldsymbol{\theta})\mathbf{V}_n(\boldsymbol{\theta})$ and $\mathbf{Q}(\boldsymbol{\theta})$ is a positively definite matrix of order $(p+1) \times (p+1)$. Moreover, we assume that

$$(2.1) \qquad \frac{1}{n}\sum_{i=1}^n \|\mathbf{v}_i(\boldsymbol{\theta})\|^4 = \mathcal{O}(1) \quad \text{as} \ \ n \to \infty$$

where $\mathbf{v}'_i(\boldsymbol{\theta})$ is the $i$-th row of $\mathbf{V}_n(\boldsymbol{\theta})$, $i = 1, \ldots, n$.

The motivation for and the validity of conditions **(A.1)**–**(A.6)** are discussed by Jurečková and Procházka [10]; they are in correspondence with the practical problems studied by these authors. Conditions **(A.3)** and **(A.6)** are also assumed by Mukherjee [18] with the difference that he replaces (2.1) with $\max_{1\leq i\leq n} \|\mathbf{v}'_i(\boldsymbol{\theta})\| = o(n^{\frac{1}{2}})$. Mukherjee, considering the time series model, assumes neither the monotonicity of $g$ in the components of $\boldsymbol{\theta}$ nor the positivity of $\mathbf{x}$, but other conditions suitable for the AR model; these conditions are also interpreted by Koul [16].

The regression $\alpha$-quantile $\widehat{\boldsymbol{\theta}}_{n\alpha}$ of model (1.1) is the minimizer of

$$(2.2) \qquad \sum_{i=1}^n \rho_\alpha(Y_i - g(\mathbf{x}_i, \mathbf{t})) = \min$$

with respect to $\mathbf{t} \in \mathbb{R}^{p+1}$, where

$$(2.3) \qquad \rho_\alpha(z) = |z|\{\alpha I[z>0] + (1-\alpha)I[z<0]\}, \ z \in \mathbb{R}^1.$$

Regarding (2.2), $\widehat{\boldsymbol{\theta}}_{n\alpha}$ can be also defined as a component $\mathbf{t}$ of the solution $(\mathbf{t}, \mathbf{r}^+, \mathbf{r}^-) \in \mathbb{R}^{p+1} \times \mathbb{R}^n_+ \times \mathbb{R}^n_+$ of the minimization

$$\alpha \sum_{i=1}^n r_i^+ + (1-\alpha)\sum_{i=1}^n r_i^- := \min$$

subject to $\quad r_i^+ - r_i^- = Y_i - g(\mathbf{x}_i, \mathbf{t}), \quad i=1,\ldots,n.$

We can also write

$$(2.4) \quad \begin{aligned} \sum_{i=1}^n \rho_\alpha(Y_i - g(\mathbf{x}_i, \widehat{\boldsymbol{\theta}}_{n\alpha})) &= \sum_{i=1}^n (Y_i - g(\mathbf{x}_i, \widehat{\boldsymbol{\theta}}_{n\alpha}))\left(I[Y_i \geq g(\mathbf{x}_i, \widehat{\boldsymbol{\theta}}_{n\alpha})]\right.\\ &\left.-(1-\alpha)\right) = \sum_{i=1}^n \left[E_{i\alpha} - [g(\mathbf{x}_i, \boldsymbol{\theta}+\mathbf{T}_n) - g(\mathbf{x}_i, \boldsymbol{\theta})]\right]\\ &\cdot \left[I[E_{i\alpha} \geq g(\mathbf{x}_i, \boldsymbol{\theta}+\mathbf{T}_n) - g(\mathbf{x}_i, \boldsymbol{\theta})] - (1-\alpha)\right] = \min, \end{aligned}$$

where

(2.5) $\quad E_{i\alpha} = e_i - F^{-1}(\alpha), \quad i = 1, \ldots, n \quad$ and

$$\mathbf{T}_n = \widehat{\boldsymbol{\theta}}_{n\alpha} - \boldsymbol{\theta}_\alpha \quad \text{and} \quad \boldsymbol{\theta}_\alpha = \boldsymbol{\theta} + \mathbf{e}_1 F^{-1}(\alpha), \quad \mathbf{e}_1 = (1, 0, \ldots, 0)'.$$

By El-Attar, Vidyasagar, and Dutta [2], if $\widehat{\boldsymbol{\theta}}_{n\alpha}$ minimizes (2.2), then there exists a vector $\mathbf{a}_n(\alpha) \in [0,1]^n$, $0 < \alpha < 1$, such that

(2.6) $\quad a_{ni}(\alpha) = \begin{cases} 1 & \text{if} \quad Y_i > g(\mathbf{x}_i, \widehat{\boldsymbol{\theta}}_{n\alpha}) \\ 0 & \text{if} \quad Y_i < g(\mathbf{x}_i, \widehat{\boldsymbol{\theta}}_{n\alpha}), \quad i = 1, \ldots, n, \end{cases}$

(2.7) $\quad \sum_{i=1}^n v_{ij}(\widehat{\boldsymbol{\theta}}_{n\alpha})[a_{ni}(\alpha) - (1-\alpha)] = 0, \quad j = 0, 1, \ldots, p,$

(2.8) $\quad \frac{1}{n}\sum_{i=1}^n (Y_i - g(\mathbf{x}_i, \widehat{\boldsymbol{\theta}}_{n\alpha}))[\hat{a}_{ni}(\alpha) - (1-\alpha)]$

$$= \frac{1}{n}\sum_{i=1}^n \rho_\alpha(Y_i - g(\mathbf{x}_i, \widehat{\boldsymbol{\theta}}_{n\alpha})).$$

**Remark 2.1.** In fact, (2.6) follows from (2.8).

Hence, we can define the regression rank scores $\hat{a}_{n1}(\alpha), \ldots, \hat{a}_{nn}(\alpha)$ as one of the vectors satisfying (2.6)–(2.8); because the set $\mathcal{A}$ of such vectors is convex, we can define $\widehat{\mathbf{a}}_n(\alpha)$ as $\mathbf{a} \in \mathcal{A}$ maximizing $\sum_{i=1}^n Y_i a_i$. Unlike Mukherjee [18], we profit from the equation (2.7). Notice that (2.7) among others implies

(2.9) $\quad \sum_{i=1}^n \hat{a}_{ni}(\alpha) = n(1-\alpha), \quad 0 < \alpha < 1,$

hence, by the continuity,

(2.10) $\quad \hat{a}_{ni}(0) = 1, \quad \hat{a}_{ni}(1) = 0, \quad i = 1, \ldots, n.$

Regarding (2.4), we can rewrite (2.6)–(2.8) in the form

(2.11) $\quad \hat{a}_{ni}(\alpha) = \begin{cases} 1 & \text{if} \quad E_{i\alpha} > g(\mathbf{x}_i, \boldsymbol{\theta} + \mathbf{T}_n) - g(\mathbf{x}_i, \boldsymbol{\theta}) \\ 0 & \text{if} \quad E_{i\alpha} < g(\mathbf{x}_i, \boldsymbol{\theta} + \mathbf{T}_n) - g(\mathbf{x}_i, \boldsymbol{\theta}), \quad i = 1, \ldots, n, \end{cases}$

(2.12) $\quad \sum_{i=1}^n v_{ij}(\widehat{\boldsymbol{\theta}}_{n\alpha})[a_{ni}(\alpha) - (1-\alpha)] = 0, \quad j = 0, 1, \ldots, p,$

$$\frac{1}{n}\sum_{i=1}^n \Big(E_{i\alpha} - [g(\mathbf{x}_i, \boldsymbol{\theta} + \mathbf{T}_n) - g(\mathbf{x}_i, \boldsymbol{\theta})]\Big)[\hat{a}_{ni}(\alpha) - (1-\alpha)]$$

(2.13) $\quad = \frac{1}{n}\sum_{i=1}^n \rho_\alpha\Big(E_{i\alpha} - [g(\mathbf{x}_i, \boldsymbol{\theta} + \mathbf{T}_n) - g(\mathbf{x}_i, \boldsymbol{\theta})]\Big).$

Jurečková and Procházka [10] proved that, under conditions **(A.1)**–**(A.6)**,

(2.14) $\quad \widehat{\boldsymbol{\theta}}_{n\alpha} - \boldsymbol{\theta}_\alpha = \mathcal{O}_p(n^{-1/2}) \;\text{ as }\; n \to \infty,$

(2.15) $\quad n^{1/2}(\widehat{\boldsymbol{\theta}}_{n\alpha} - \boldsymbol{\theta}_\alpha)$

$$= \frac{1}{n^{1/2}f(F^{-1}(\alpha))} \sum_{i=1}^n \mathbf{v}_i(\boldsymbol{\theta})\psi_\alpha(E_{i\alpha}) + o_p(1) \;\text{ as }\; n \to \infty$$

uniformly in $\varepsilon \leq \alpha \leq 1 - \varepsilon$, $\forall \varepsilon \in (0, \frac{1}{2})$, where

(2.16) $$\psi_\alpha(u) = I[u \geq 0] - (1 - \alpha), \ u \in \mathbb{R}.$$

Expanding $g(\mathbf{x}_i, \widehat{\boldsymbol{\theta}}_{n\alpha}) - g(\mathbf{x}_i, \boldsymbol{\theta}_\alpha)$ around $\widehat{\boldsymbol{\theta}}_{n\alpha}$, we obtain for $i = 1, \ldots, n$

$$g(\mathbf{x}_i, \widehat{\boldsymbol{\theta}}_{n\alpha}) - g(\mathbf{x}_i, \boldsymbol{\theta}_\alpha) = \widehat{\theta}_{0,n\alpha} - \theta_0 - F^{-1}(\alpha) + \tilde{g}(\mathbf{x}_i, \widehat{\boldsymbol{\theta}}^*_{n\alpha}) - \tilde{g}(\mathbf{x}_i, \boldsymbol{\theta}^*)$$

(2.17) $$= \widehat{\theta}_{0,n\alpha} - \theta_0 - F^{-1}(\alpha) + \sum_{j=1}^{p} v_{ij}(\widehat{\theta}_{j,n\alpha})(\widehat{\theta}_{j,n\alpha} - \theta_j)$$

$$-(\widehat{\boldsymbol{\theta}}^*_{n\alpha} - \boldsymbol{\theta}^*)' \mathbf{B}^{(i)}_*(\widetilde{\boldsymbol{\theta}}_n)(\widehat{\boldsymbol{\theta}}^*_{n\alpha} - \boldsymbol{\theta}^*)/2$$

with $\widetilde{\boldsymbol{\theta}}_n$ between $\boldsymbol{\theta}^*$ and $\widehat{\boldsymbol{\theta}}^*_{n\alpha}$, and with

$$\mathbf{B}^{(i)}_*(\widetilde{\boldsymbol{\theta}}_n) = \left[\frac{\partial^2 g(\mathbf{x}_i, \boldsymbol{\theta})}{\partial \theta_j \partial \theta_k}\right]^p_{j,k=1}\bigg|_{\boldsymbol{\theta}=\widetilde{\boldsymbol{\theta}}}.$$

Inserting (2.17) in (2.13) and regarding (2.11), (2.14) and condition **(A.1)**, we obtain

$$\frac{1}{n}\sum_{i=1}^{n}\Big(E_{i\alpha} - [g(\mathbf{x}_i, \boldsymbol{\theta} + \mathbf{T}_n) - g(\mathbf{x}_i, \boldsymbol{\theta})]\Big)[\hat{a}_i(\alpha) - (1-\alpha)]$$

(2.18) $$= \frac{1}{n}\sum_{i=1}^{n} E_{i\alpha}[\hat{a}_i(\alpha) - (1-\alpha)] + \alpha(1-\alpha)\cdot\mathcal{O}_p(n^{-1})$$

uniformly in $\varepsilon \leq \alpha \leq 1-\varepsilon$, $\varepsilon \in (0, \frac{1}{2})$. On the other hand, it follows from Lemma 3.5 in Jurečková and Procházka [10] and its proof that

$$\frac{1}{n}\sum_{i=1}^{n}\rho_\alpha\Big(e_i - [g(\mathbf{x}_i, \boldsymbol{\theta} + \mathbf{T}_n) - g(\mathbf{x}_i, \boldsymbol{\theta})]\Big)$$

$$= \frac{1}{n}\sum_{i=1}^{n}\rho_\alpha(E_{i\alpha}) + \alpha(1-\alpha)\cdot\mathcal{O}_p(n^{-1})$$

(2.19) $$= \frac{1}{n}\sum_{i=1}^{n} E_{i\alpha}\Big(I[e_i \geq F^{-1}(\alpha)] - (1-\alpha)\Big) + \alpha(1-\alpha)\cdot\mathcal{O}_p(n^{-1})$$

uniformly in $\varepsilon \leq \alpha \leq 1 - \varepsilon$, $\varepsilon \in (0, \frac{1}{2})$. Combining (2.11), (2.18) and (2.19) entails

(2.20) $$\frac{1}{n}\sum_{i=1}^{n} E_{i\alpha}\Big(\hat{a}_i(\alpha) - I[e_i \geq F^{-1}(\alpha)]\Big) = \alpha(1-\alpha)\cdot\mathcal{O}_p(n^{-1})$$

uniformly in $\varepsilon \leq \alpha \leq 1 - \varepsilon$, $\varepsilon \in (0, \frac{1}{2})$. This, in turn, further implies

(2.21) $$\frac{1}{n}\sum_{i=1}^{n} E_{i\alpha}[\hat{a}_i(\alpha) - (1-\alpha)] = \frac{1}{n}\sum_{i=1}^{n} E_{i\alpha}[a_n^*(R_{ni}, \alpha) - (1-\alpha)] + \alpha(1-\alpha)\cdot\mathcal{O}(n^{-1})$$

uniformly in $\varepsilon \leq \alpha \leq 1 - \varepsilon$, $\varepsilon \in (0, \frac{1}{2})$, where $a_n^*(R_{ni}, \alpha)$ are Hájek's rank scores,

(2.22) $$a_n^*(R_{ni}, \alpha) = \begin{cases} 0, & \text{if } \frac{R_{ni}}{n} < \alpha, \\ R_{ni} - \alpha & \text{if } \frac{R_{ni}-1}{n} \leq \alpha < \frac{R_{ni}}{n} \\ 1 & \text{if } \alpha < \frac{R_{ni}-1}{n} \end{cases}$$

and $R_{ni}$ is the rank of $e_i$, $i = 1, \ldots, n$ (see, e.g., Hájek and Šidák [5], Section V.3.5).

Consider a triangular array $\{\mathbf{z}_{n1}, \ldots, \mathbf{z}_{nn}\}$ of vectors from $\mathbb{R}^r$ such that

**(Z.1)** $\sum_{i=1}^n z_{ni} = 0$ and $\frac{1}{n}\sum_{i=1}^n \mathbf{z}_{ni}\mathbf{z}'_{ni} \to \mathbf{C}$ as $n \to \infty$ where $\mathbf{C}$ is a positive $r \times r$ matrix,

**(Z.2)** $\max_{1 \le i \le n} \|\mathbf{z}_{ni}\| = o(n^{1/2})$,

and the nonlinear rank statistic process

$$(2.23) \qquad \mathbf{Z}_n(\alpha) = \frac{1}{n}\sum_{i=1}^n \mathbf{z}_{ni}[\hat{a}_{ni}(\alpha) - (1-\alpha)], \ 0 < \alpha < 1.$$

Let $\varphi : (0,1) \mapsto \mathbb{R}$ be a monotone function. Fix $\varepsilon \in (0, \frac{1}{2})$ and put

$$(2.24) \qquad \varphi_\varepsilon(u) = \begin{cases} \varphi(\varepsilon) & \text{if} \quad 0 \le u < \varepsilon \\ \varphi(u) & \text{if} \quad \varepsilon \le u \le 1-\varepsilon \\ \varphi(1-\varepsilon) & \text{if} \quad 1-\varepsilon < u \le 1 \end{cases}$$

and define the scores

$$(2.25) \qquad \hat{b}_{ni} = \int_0^1 \hat{a}_{ni}(\alpha)d\varphi_\varepsilon(\alpha), \ i = 1, \ldots, n.$$

For example, $\varphi(u) = u - \frac{1}{2}$ (Wilcoxon score function) or

$$\varphi(u) = \begin{cases} -1 & \text{if} \quad 0 \le u < \frac{1}{2} \\ 0 & \text{if} \quad u = \frac{1}{2} \\ 1 & \text{if} \quad \frac{1}{2} < u \le 1 \end{cases}$$

(median score function). The vector of nonlinear rank statistics

$$(2.26) \qquad \mathbf{S}_n = \sum_{i=1}^n \mathbf{z}_{ni}\hat{b}_{ni}$$

will serve us as a basis for a construction of tests.

**Remark 2.2.** Gutenbrunner et al. [4] and Hallin and Jurečková [6] considered the scores of type (2.25) with $\varphi_\varepsilon$ replaced by a nondecreasing, square-integrable score function $\varphi : (0,1) \mapsto \mathbb{R}$, satisfying conditions of Chernoff–Savage type, including the normal scores. However, they were not able to construct the tests with nuisance linear regression for $f$ with heavy tails. Jurečková [9] admitted heavy-tailed distributions, but her scores were truncated to $(\varepsilon_n, 1 - \varepsilon_n)$ with $\varepsilon_n \downarrow 0$. Under the conditions **(A.1)**–**(A.3)** and using the methods of Jurečková and Procházka [10], we are able to guarantee the uniformity in (2.15) only on $[\varepsilon, 1-\varepsilon]$ with a fixed $\varepsilon \in (0, \frac{1}{2})$. On the other hand, we do not restrict the tails of the distribution.

A possible extension of the subinterval $[\varepsilon, 1-\varepsilon]$ to $[\varepsilon_n, 1-\varepsilon_n]$ or to $(0,1)$ will be an object of a further study.

## 3. Properties of nonlinear rank statistics

We shall first prove the lemma.

**Lemma 3.1.** *Let* $\{\mathbf{z}_{n1}, \ldots, \mathbf{z}_{nn}\}$ *be a triangular array of vectors from* $\mathbb{R}^r$ *satisfying* **(Z.1)** *and* **(Z.2)**, *let*

$$\mathbf{Z}_n = \begin{bmatrix} \mathbf{z}'_{n1} \\ \cdots \\ \mathbf{z}'_{nn} \end{bmatrix}.$$

*Then, under conditions* **(A.1)**–**(A.3)**,

$$\sup_{\varepsilon \leq \alpha \leq 1-\varepsilon} \left\{ n^{-1/2} \Big\| \sum_{i=1}^n \Big[ \mathbf{z}_{ni} \psi_\alpha(E_{i\alpha} - [g(\mathbf{x}_i, \boldsymbol{\theta} + \mathbf{T}_n) - g(\mathbf{x}_i, \boldsymbol{\theta})]) \right.$$

(3.1)
$$\left. - (\mathbf{z}_{ni} - \hat{\mathbf{z}}_{ni}) \psi_\alpha(E_{i\alpha}) \Big] \Big\| \right\} \xrightarrow{p} \mathbf{0}$$

*as* $n \to \infty$, *for any fixed* $\varepsilon \in (0, \frac{1}{2})$, *where* $\hat{\mathbf{z}}'_{ni}$ *is the i-th row of the projection* $\widehat{\mathbf{Z}}_n$ *of* $\mathbf{Z}_n$ *in the space spanned by the columns of matrix* $\mathbf{V}_n(\boldsymbol{\theta})$, *i.e.*

(3.2)
$$\widehat{\mathbf{Z}}_n = \widehat{\mathbf{H}}_n(\boldsymbol{\theta}) \mathbf{Z}_n,$$

$$\widehat{\mathbf{H}}_n(\boldsymbol{\theta}) = \mathbf{V}_n(\boldsymbol{\theta}) \left[ \mathbf{V}'_n(\boldsymbol{\theta}) \mathbf{V}_n(\boldsymbol{\theta})' \right]^{-1} \mathbf{V}'_n(\boldsymbol{\theta}).$$

*Proof.* It follows from Lemma 3.5 in Jurečková and Procházka [10] that

(3.3)
$$\frac{1}{n} \sum_{i=1}^n [\rho_\alpha(E_{i\alpha} - g(\mathbf{x}_i, \boldsymbol{\theta} + \mathbf{t})) - \rho_\alpha(E_{i\alpha})] = \frac{1}{2} f(F^{-1}(\alpha)) \mathbf{t}' \mathbf{Q}_n \mathbf{t}$$

(3.4)
$$-\mathbf{t}' \frac{1}{n} \sum_{i=1}^n \mathbf{v}_i(\boldsymbol{\theta}) \psi_\alpha(E_{i\alpha}) + \mathcal{O}_p\left(n^{-1/2} \|\mathbf{t}\|^{3/2} + \|\mathbf{t}\|^2\right)$$

uniformly in $\varepsilon \leq \alpha \leq 1 - \varepsilon$ and in $\|\mathbf{t}\| \leq r_n$, for every sequence $\{r_n\}$ of positive numbers tending to 0. Now, extend the model (1.1) in the following way:

(3.5)
$$Y_i = g(\mathbf{x}_i, \boldsymbol{\theta}) + \mathbf{z}'_{ni} \boldsymbol{\vartheta} + e_i, \quad i = 1, \ldots, n, \quad \boldsymbol{\vartheta} \in \mathbb{R}^r.$$

Then the conditions **(A.1)**–**(A.3)** are satisfied even for extended model (3.5), replacing $\boldsymbol{\theta}$ by $(\boldsymbol{\theta}', \boldsymbol{\vartheta}')'$. The function $\rho_\alpha$ is absolutely continuous and convex; taking the right derivative of (3.3) with respect to last $r$ coordinates of $\mathbf{t}$ (evaluated when the last $r$ coordinates of $\mathbf{t}$ are zero), we obtain

$$\sup_{\varepsilon \leq \alpha \leq 1-\varepsilon} \left\{ n^{-1/2} \Big\| \sum_{i=1}^n \Big[ \mathbf{z}_{ni} \psi_\alpha(E_{i\alpha} - [g(\mathbf{x}_i, \boldsymbol{\theta} + \mathbf{t}) - g(\mathbf{x}_i, \boldsymbol{\theta})]) \right.$$

(3.6)
$$\left. - (\mathbf{z}_{ni} - \hat{\mathbf{z}}_{ni}) \psi_\alpha(E_{i\alpha}) \Big] \Big\| \right\} = o_p(1)$$

uniformly in $\mathbf{t} \in \mathbb{R}^{p+1}$, $\|\mathbf{t}\| \leq r_n$, for every sequence $\{r_n\}$ of positive numbers tending to 0. Inserting $\mathbf{t} \mapsto \widehat{\boldsymbol{\theta}}_{n\alpha} - \boldsymbol{\theta}_\alpha$ into (3.6), we arrive at (3.1). □

The following corollary approximates the regression rank scores by an empirical process.

**Corollary 3.1.** *Under the conditions of Lemma 3.1,*

(3.7)
$$\sup_{\varepsilon \leq \alpha \leq 1-\varepsilon} \left\{ n^{-1/2} \Big\| \sum_{i=1}^n \Big[ \mathbf{z}_{ni} \hat{a}_{ni}(\alpha) - (\mathbf{z}_{ni} - \hat{\mathbf{z}}_{ni}) I[e_i \geq F^{-1}(\alpha)] \Big] \Big\| \right\} \xrightarrow{p} 0$$

*for any fixed* $\varepsilon \in (0, \frac{1}{2})$.

*Proof.* Notice that

$$\sup_{\varepsilon \leq \alpha \leq 1-\varepsilon} \left\{ n^{-1/2} \Big\| \sum_{i=1}^{n} \mathbf{z}_{ni} I[Y_i = g(\mathbf{x}_i, \widehat{\boldsymbol{\theta}}_{n\alpha})] \Big\| \right\} \xrightarrow{p} 0;$$

hence, regarding (2.6), (3.7) follows from Lemma 3.1. $\square$

**Corollary 3.2.** *Under the conditions of Lemma 3.1,*

$$(3.8) \quad \sup_{\varepsilon \leq \alpha \leq 1-\varepsilon} \left\{ n^{-1/2} \Big\| \sum_{i=1}^{n} \Big[ \mathbf{z}_{ni} \hat{a}_{ni}(\alpha) - (\mathbf{z}_{ni} - \hat{\mathbf{z}}_{ni}) a_n^*(R_{ni}, \alpha) \Big] \Big\| \right\} \xrightarrow{p} 0$$

*as $n \to \infty$, for any fixed $\varepsilon \in (0, \frac{1}{2})$, where $a_n^*(R_{ni}, \alpha)$ are Hájek's scores defined in (2.22).*

Let $\varphi : (0,1) \mapsto \mathbb{R}$ be a monotone function, fix $\varepsilon \in (0, \frac{1}{2})$ and consider the scores (2.25) and the nonlinear rank statistics (2.26). Then Corollaries 3.1, 3.2 imply

**Corollary 3.3.** *Under the conditions of Lemma 3.1, for any fixed $\varepsilon \in (0, \frac{1}{2})$,*

$$n^{-1/2} \Big\| \sum_{i=1}^{n} \Big[ \mathbf{z}_{ni} \hat{b}_{ni} - (\mathbf{z}_{ni} - \hat{\mathbf{z}}_{ni}) \int_0^1 I[e_i \geq F^{-1}(\alpha)] d\varphi_\varepsilon(\alpha) \Big] \Big\| = o_p(1),$$

$$(3.9) \quad n^{-1/2} \Big\| \sum_{i=1}^{n} \Big[ \mathbf{z}_{ni} \hat{b}_{ni} - (\mathbf{z}_{ni} - \hat{\mathbf{z}}_{ni}) \hat{b}_{ni}^* \Big] \Big\| = o_p(1),$$

*where*

$$\hat{b}_{ni}^* = \int_0^1 a_n^*(R_{ni}, \alpha) d\varphi_\varepsilon(\alpha), \quad i = 1, \ldots, n.$$

*Hence, under the model (1.1), the nonlinear rank statistic (2.26) is asymptotically equivalent to the linear rank statistic $n^{-1/2} \sum_{i=1}^{n} (\mathbf{z}_{ni} - \hat{\mathbf{z}}_{ni}) \hat{b}_{ni}^*$ pertaining to the Hájek rank scores.*

## 4. Application: Tests of linear regression in nonlinear regression model with unknown parameters

Convenient properties of the nonlinear regression rank scores, proved in the former sections, lead to an idea of their possible application in testing the significance of a linear regression in the presence of a nonlinear regression with nuisance parameters, or in testing other hypotheses with nuisance parameters of nonlinear regression. For instance, we can compare two sets of observations affected by a nonlinear regression with unknown parameters.

Let us illustrate possible applications on the nonlinear regression model

$$(4.1) \quad Y_i = g(\mathbf{x}_{ni}, \boldsymbol{\theta}) + \mathbf{z}_{ni}' \boldsymbol{\beta} + e_i, \quad i = 1, \ldots, n$$

where $\boldsymbol{\beta} \in \mathbb{R}^r$ is an unknown parameter and $\mathbf{z}_{ni}$, $i = 1, \ldots, n$, are known (or observable) regressors. Denote

$$\mathbf{Z}_n = \begin{bmatrix} \mathbf{z}_{n1}' \\ \ldots \\ \mathbf{z}_{nn}' \end{bmatrix}$$

the matrix of order $n \times r$ and assume that it has rank $r$. The problem is that of testing the hypothesis

$$\mathbf{H}_0 : \boldsymbol{\beta} = \mathbf{0}$$

with $g(\cdot, \cdot)$ of known shape and with known $\mathbf{x}_i$, $i = 1, \ldots, n$, but with $\boldsymbol{\theta}$ and $F$ unspecified; we only assume that conditions **(A.1)**–**(A.3)** and **(Z.1)**–**(Z.2)** are satisfied.

As an example consider the situation where the nonlinear regression describes the concentration of the drug teofylin in the human blood, measured at times $t_1, \ldots, t_k$ after the application. In such situation, the regression function is typically a mixture of exponentials with unknown parameters (see Jurečková and Procházka [10] and Schindler [19]). Treating two groups of patients (boys and girls), we want to compare their reactions to the drug. This leads to a two-sample problem with a nuisance nonlinear regression.

Let $\hat{a}_{n1}(\alpha), \ldots, \hat{a}_{nn}(\alpha)$, $0 \leq \alpha \leq 1$ denote the regression rank scores corresponding to the submodel under $\mathbf{H}_0$,

(4.2) $$Y_i = g(\mathbf{x}_i, \boldsymbol{\theta}) + e_i, \ i = 1, \ldots, n.$$

Let $\varphi : (0, 1) \mapsto \mathbb{R}$ be a nondecreasing, bounded score function such that $\varphi(1 - u) = -\varphi(u)$, $0 < u < 1$. Fix $\varepsilon \in (0, \frac{1}{2})$, calculate the scores $\hat{b}_{ni}$, $i = 1, \ldots, n$, defined in (2.25), and the test criterion

(4.3) $$\mathcal{T}_n = (A(\varphi_\varepsilon))^{-2} \mathbf{S}'_n \mathbf{D}_n^{-1} \mathbf{S}_n, \quad \mathbf{S}_n = n^{-1/2} \sum_{i=1}^n \mathbf{z}_{ni} \hat{b}_{ni}$$

where

$$(A(\varphi_\varepsilon))^2 = \int_0^1 (\varphi_\varepsilon(u) - \bar{\varphi}_\varepsilon)^2 du, \quad \bar{\varphi}_\varepsilon = \int_0^1 \varphi_\varepsilon(u) du$$

and

$$\mathbf{D}_n = \frac{1}{n} \left( \mathbf{Z}_n - \widehat{\mathbf{Z}}_n \right)' \left( \mathbf{Z}_n - \widehat{\mathbf{Z}}_n \right),$$

and $\widehat{\mathbf{Z}}_n$ is the projection of $\mathbf{Z}_n$ in the space spanned by the columns of matrix $\mathbf{V}_n(\widehat{\boldsymbol{\theta}}_n)$ with $\widehat{\boldsymbol{\theta}}_n$ being the nonlinear regression quantile of model (4.2) with $\alpha = \frac{1}{2}$, i.e.

$$\widehat{\mathbf{Z}}_n = \widehat{\mathbf{H}}_n(\widehat{\boldsymbol{\theta}}_n) \mathbf{Z}_n,$$
$$\widehat{\mathbf{H}}_n(\widehat{\boldsymbol{\theta}}_n) = \mathbf{V}_n(\widehat{\boldsymbol{\theta}}_n) \left( \mathbf{V}'_n(\widehat{\boldsymbol{\theta}}_n) \mathbf{V}'_n(\widehat{\boldsymbol{\theta}}_n) \right)^{-1} \mathbf{V}'_n(\widehat{\boldsymbol{\theta}}_n).$$

The test is based on the asymptotic distribution of $\mathcal{T}_n$ under $\mathbf{H}_0$. Regarding that

$$\mathbf{V}_n(\widehat{\boldsymbol{\theta}}_n) - \mathbf{V}_n(\boldsymbol{\theta}) \xrightarrow{p} \mathbf{0} \ \text{as} \ n \to \infty,$$

the asymptotic distribution of $\mathcal{T}_n$ under $\mathbf{H}_0$, due to the Corollary 3.1, in turn coincides with the asymptotic distribution of the sequence

$$(A(\varphi_\varepsilon))^{-2} \mathbf{S}_n^{*\prime} \mathbf{D}_n^{-1} \mathbf{S}_n^*, \quad \mathbf{S}_n^* = n^{-1/2} \sum_{i=1}^n (\mathbf{z}_{ni} - \hat{\mathbf{z}}_{ni}) \varphi_\varepsilon \left( \frac{R_{ni}}{n+1} \right),$$

where $R_{ni}$ is the rank of $e_i$, $i = 1, \ldots, n$. Hence, the test of $\mathbf{H}_0$ in the presence of nuisance nonlinear regression would coincide with the ordinary rank test under $\boldsymbol{\theta}$ known, with score generating function $\varphi_\varepsilon$ and with the coefficients $\mathbf{z}_{ni} - \hat{\mathbf{z}}_{ni}$, $i = 1, \ldots, n$. The asymptotic null distribution of $\mathcal{T}_n$ under $\mathbf{H}_0$ is the central $\chi^2$ with $r$ degrees of freedom. The asymptotic distribution of $\mathcal{T}_n$ under the local alternative

$$\mathbf{H}_n: \ \boldsymbol{\beta}_n = n^{-1/2}\boldsymbol{\beta}_0 \quad (\boldsymbol{\beta}_0 \in \mathbb{R}^r \ \text{fixed})$$

is the noncentral $\chi^2$ with $r$ degrees of freedom and with noncentrality parameter

$$(A(\varphi_\varepsilon))^{-2} \ \boldsymbol{\beta}_0'\mathbf{D}\boldsymbol{\beta}_0 \left[\int_0^1 \varphi_\varepsilon(u) df(F^{-1}(u))\right]^2, \quad \mathbf{D} = \lim_{n\to\infty} \mathbf{D}_n.$$

If $\beta \in \mathbb{R}^1$, i.e. $r = 1$, we can test $\mathbf{H}_0: \ \beta = 0$ also against the one-sided alternative $\mathbf{H}_1: \ \beta > 0$, what is most interesting in the two-sample model. Then the test criterion simplifies to

$$(4.4) \qquad \mathcal{T}_n^* = (A(\varphi_\varepsilon))^{-1} \Big[\sum_{i=1}^n (z_{ni} - \hat{z}_{ni})^2\Big]^{-1/2} \sum_{i=1}^n z_{ni}\hat{b}_{ni}$$

and rejects $\mathbf{H}_0$ in favor of $\mathbf{H}_1$ on the asymptotic significance level $\tau$ provided $\mathcal{T}_n^*$ exceeds the $(1-\tau)$-quantile of the standard normal distribution, i.e. $\mathcal{T}_n^* \geq \Phi^{-1}(1-\tau)$. The performance of such test on the real data is illustrated by Schindler [19], who has also elaborated a suitable computation algorithm. Neither the Wilcoxon nor the median tests indicate a difference in the dynamism of spreading the teofylin between two groups of patients.

A general class of tests using the regression rank scores, including their numerical behavior and the algorithms, is a subject of Schindler's PhD Thesis (Schindler [20]). Besides the multiple comparisons, we can also test for the independence of two random variables, affected by a nonlinear regression with unknown parameters, or test for the significance of an outward trend in a similar situation, and we have various other applications.

**Acknowledgment.** I would like to thank Pranab Kumar Sen for all our cooperation, which has dated since 1980. I really learned much from him, and I am still learning.